\documentclass{amsart}
\usepackage{amsfonts, amssymb, amsthm}
\theoremstyle{plain}
\newtheorem{theorem}{Theorem}

\newtheorem{vardefinition}[theorem]{Definition}
\newtheorem{lemma}[theorem]{Lemma}
\newtheorem{proposition}[theorem]{Proposition}

\newtheorem{varremark}[theorem]{Remark}
\newtheorem{varexample}[theorem]{Example}
\newtheorem{varexercise}[theorem]{Exercise}

\newenvironment{definition}{\begin{vardefinition}
\begin{normalfont}}{\end{normalfont}
\end{vardefinition}}

\newenvironment{example}{\begin{varexample}
\begin{normalfont}}{\end{normalfont}
\end{varexample}}

\DeclareMathOperator{\Ass}{Ass}
\DeclareMathOperator{\Spec}{Spec}

\newcommand{\m}{\mathfrak{m}}

\numberwithin{theorem}{section}

\newcommand{\suchthat}{\; | \;}

\newcommand{\N}{{\mathbb N}}

\title{Semilocal Formal Fibers
of Principal Prime Ideals}

\thanks{This research was supported by National Science Foundation grant DMS-0353634. The authors would like to acknowledge S. Loepp for encouragement and many useful suggestions.}

\author{John Chatlos}
\address{J. Chatlos: Williams College}
\email{John.G.Chatlos@williams.edu}

\author{Brian Simanek}
\address{B. Simanek: California Institute of Technology}
\email{bsimanek@caltech.edu}

\author{Nathaniel G. Watson}
\address{N. G. Watson: The University of California, Berkeley}
\email{watson@math.berkeley.edu}

\author{Sherry X. Wu }
\address{S. X. Wu.: Massachusetts Institute of Technology}
\email{sxw2@mit.edu}

\begin{document}


\begin{abstract}
Let $(T,\m)$ be a complete local (Notherian) ring,
$C$ a finite set of pairwise incomparable nonmaximal prime ideals of $T$, and $p \in T$ a nonzero element.
We provide necessary and sufficient conditions for $T$ to be the completion of an integral domain $A$
containing the prime ideal $pA$ whose formal fiber is semilocal with maximal ideals the elements of $C$.
\end{abstract}

\maketitle

\section{Introduction}
One way to better understand the relationship between a commutative local ring and its completion is to examine the formal
fibers of the ring. Given  a local ring $A$ with maximal ideal $\m$ and $\m$-adic completion $\widehat{A}$, the formal
fiber of a prime ideal $P \in \Spec A$ is defined to be $\Spec(\widehat{A} \otimes_A k(P))$, where $k(P):=A_P/PA_P$. Because
there is a one-to-one correspondence between the elements in the formal fiber of  $P$ and the
prime ideals in the inverse image of $P$ under the map from $\Spec \widehat{A}$  to $\Spec A $ given by $Q \rightarrow Q \cap A$, we can think
of $Q \in \Spec \widehat{A} $ as being in the formal fiber of $P$ if and only if $Q\cap A=P$.

One fruitful way of researching formal fibers has been, instead of directly computing the formal fibers of rings, to investigate ``inverse" formal fiber questions---that is, given a complete local ring $T,$ when does there exist a local ring $A$ such that $\widehat{A}=T$ and both $A$ and the formal fibers of prime ideals in $A$ meet certain prespecified conditions? One important result of this type is due to P.\ Charters and S.\ Loepp, who show in \cite{SLGFF} that given a complete local ring $T$ with maximal ideal $\m$ and $G \subset \Spec T $ where $G$ is a finite set of prime ideals which are pairwise incomparable by inclusion, there exists a local domain $A$ such that $\widehat{A}=T$ and the formal fiber of the zero ideal of $A$ is semilocal with maximal ideals exactly the elements of $G$
if and only if certain relatively weak conditions are satisfied.

In this paper we address a similar question: what are the
necessary and sufficient conditions for $T$ to be the completion of a local domain $A$ possessing a principal prime ideal with a specified semilocal formal fiber?

Partial results on this subject were achieved by A.\ Dundon, D.\ Jensen, S. \ Loepp, J.\ Provine, and J.\ Rodu in \cite{SM03}, under the constraint that the specified set $G=\{Q_1, Q_2,\ldots, Q_k\}$ of nonmaximal ideals in the formal fiber is such that $\bigcap_{i=1}^k Q_i$ contains a nonzero regular {\em prime} element $p$ of $T$. In particular, suppose this holds and with $\Pi$ denoting the prime subring of $T,$ we have that either $\Pi \cap Q_i =(0)$ for every $i$ or $\Pi \cap Q_i =p\Pi$ for every $i.$ In \cite{SM03} it
is shown that there exists a local
domain $A$ such that $\widehat{A}=T$, $p \in A$, $pA \in \Spec A$ and the formal fiber of $pA$ is semilocal with maximal ideals the elements of $C$ if and only if the following conditions are satisfied:
\begin{enumerate}
\item $p \in Q_i$ for every $i$.
\item If $\dim T=1$ then ${C}={\m}$.
\item If $\dim T >1$, then $\m \notin {C}$.
\end{enumerate}
Note that the condition that $T$ contains the prime ideal $pT$ for a nonzero $p$ in $T$ implies that $T$ is a domain.

The main theorem in this paper is an improvement on the results in \cite{SM03}. We eliminate the assumption that $p$ is a prime element in $T$. Moreover, Theorem \ref{principal} in this paper provides necessary and sufficient conditions for a complete local ring to be the completion of an integral
domain containing a height one
principal prime ideal with specified semilocal formal fiber. Specifically,
let $T$ be a complete local ring with maximal ideal $\m$, $\Pi$ the prime
subring of $T$, and $C=\{Q_1, Q_2, \ldots, Q_k\}$  a finite set of nonmaximal incomparable prime ideals of $T$.  Let $p \in \bigcap_{i=1}^k Q_i$ with $p \neq 0$.  We show that there exists a local domain
$A$ with $p \in A$  such that $\widehat{A} = T$  and $pA$ is a prime ideal whose
formal fiber is semilocal with maximal
ideals the elements of $C$ if and only if the following conditions are satisfied:
\begin{enumerate}
\item $Q_i$ has height at least one for every $i \in \{1,2, \ldots, k\}$.
\item For all $P \in {\Ass T/pT}$, $ P \subseteq Q_i$ for some $i \in \{1,2, \ldots, k\}$.
\item $P \cap \Pi[p] =(0)$ for all $P \in \Ass T.$
\item $((Q_i \setminus pT) \Pi[p]) \cap \Pi[p]  = \{0\}$  for all $i \in \{1,2, \ldots, k\}.$
\end{enumerate}

The proof that the above conditions are necessary is relatively short. Therefore most of this paper is devoted to showing they are sufficient by constructing an integral domain $A$ with the desired properties. The general
strategy behind our construction, which is similar to constructions in both \cite{SLGFF} and \cite{SM03}, is to start with the prime subring of $T$ localized at its maximal ideal
and recursively build up an ascending chain of subrings maintaining some specific properties.
Our final ring $A$ will be the union of all the subrings in the chain. Most of the work in the construction goes toward insuring that $A$ simultaneously meets three conditions: the map $A \rightarrow T/J$ is onto for every ideal $J$ such that $J \nsubseteq Q_i$ for all $i \in \{1,2, \ldots, k\}$; $IT \cap A=I$ for every finitely generated ideal $I$ of $A$; and $(Q_i\setminus pT)A \cap A= \{0\}$ for all $i$. These conditions will ensure that $\widehat{A} =T$ and that $pA \in \Spec A$ has a semilocal formal fiber with maximal ideals precisely the elements of $C$.

Throughout this paper, all rings will be commutative with unity. When we say a ring is ``quasilocal" we mean that it has one maximal ideal. A ``local" ring will be a Noetherian quasilocal ring.

\section{Semilocal Formal Fibers of Principle Prime Ideals of a Domain}

Suppose we are given a complete local ring $(T,\m)$, and a finite set  $C=\{Q_1, Q_2, \ldots, Q_k\} \subseteq \Spec T$ of pairwise incomparable (that is, $Q_i \subseteq Q_j $ if and only if $Q_i = Q_j$) nonmaximal prime ideals. In this section we answer the following question.  When is it true that there is a local domain $A$ such that $\widehat{A} =T$ and  there is some principle prime $P \in \Spec A$ such that the formal fiber of $P$ is semilocal with maximal ideals $\{Q_1, Q_2, \ldots, Q_k\}$?

\begin{definition} Let $S$ be a set. Define $\Gamma(S)= \sup( |S|, \aleph_0)$.
\end{definition}
Note that clearly if $T$ and $S$ are sets, $\Gamma(S)\Gamma(T) = \sup(\Gamma(S), \Gamma(T)).$ This definition simplifies the statement of some of our lemmas.

\begin{definition}
Let $(T,\m)$ be a complete local ring and suppose we have a finite, pairwise incomparable set $C=\{Q_1, Q_2, \ldots, Q_k\} \subseteq \Spec T$.  Let $p \in \bigcap_{i=1}^{k} Q_i$ be a
nonzero regular element of $T$. Suppose that $(R,R\cap \m)$ is a quasilocal subring of $T$
containing $p$ with the following properties:
\begin{enumerate}
\item $\Gamma(R) <|T|$;
\item If $P$ is an associated prime ideal of $T$ then $R\cap P=(0)$;
\item For all $i \in \{1,2, \ldots, k\}$, $((Q_i \setminus pT) R) \cap R = \{ 0 \}$ where $(Q_i \setminus pT) R= \{ t \in T | t = qr  \mbox{ for some } q \in Q_i \setminus pT, \  r \in R\}$\footnotemark[1]

\end{enumerate}
Then we call $R$ a \emph{$pT$-complement avoiding subring of $T$}, which we shorten to
\emph{$p$ca subring}.
\end{definition}

\footnotetext[1]{Though this condition admittedly looks a little strange, if ${Q_i}^n$ is primary for each $i$ and all $n$, it is equivalent to ${Q_i}^n \cap R = p^n T \cap R$ for all $n.$}
%

To show the existence of our
local domain $A$, we construct a chain of intermediate $p$ca subrings
and then let $A$ be the union of these subrings.  For some steps of the construction we need the additional condition that $pT \cap R=pR$ for our subring $R$.  The following
lemmas show that given a $p$ca subring $R$, we can find a larger
$p$ca subring $S$ with this property. Note that these lemmas (\ref{Claim1}, \ref{Claim2}, and \ref{fixing}) are parallel to Lemmas 3.2, 3.4, and 3.5 respectively in \cite{SM03}, with the $p$in subrings of \cite{SM03} replaced with $p$ca subrings in this paper. This change necessitates that the proof of Lemma \ref{Claim1} differs substantially from the proof of Lemma 3.2 of \cite{SM03}, but the proofs of Lemma \ref{Claim2} and Lemma \ref{fixing} are essentially the same as the proofs of Lemma 3.4 and Lemma 3.5 of \cite{SM03} respectively.


\begin{lemma}
\label{Claim1}
Let $(T,\m)$ be a complete local ring and suppose we have a finite, pairwise incomparable set $C=\{Q_1, Q_2, \ldots, Q_k\} \subseteq \Spec T$.  Let $p \in \bigcap_{i=1}^{k} Q_i$ be a
nonzero regular element of $T$.  Let $(R,R \cap \m)$ be a $p$ca subring of $(T,\m)$ and let $c \in pT \cap R$.  Then there
exists a $p$ca subring $S$ of $T$ such that $R\subseteq S\subseteq T$,
$c \in pS$, and $\Gamma(S) =\Gamma(R).$
\end{lemma}

\begin{proof}
Since  $c \in pT \cap R$, $c=pu$ for some element $u$ in $T$.  We
claim that $S=R[u]_{R[u]\cap \m}$ is the desired subring.  Clearly $\Gamma(S) = \Gamma(R)$ and in particular $\Gamma(S)< |T|$.

First we consider an arbitrary $f \in R[u]$ with $f \neq 0.$ We can write $f = r_n u^n+ \cdots +r_1 u +r_0$ for some $r_0, r_1 , \ldots ,r_n \in R.$ Then $$p^n f = r_n (pu)^n + p r_{n-1} (pu)^{n-1} + \cdots+p^{n-1} r_1 (pu) + p^n r_0 = r_n c^n + p r_{n-1} c^{n-1} + \cdots+ p^{n-1} r_1 c +p^n r_0 $$ and thus we see $p^n f \in R.$

Now for any $i \in \{1,2, \ldots ,k\}$, choose an $f \in(( Q_i \setminus pT)R[u])\cap R[u].$ Since $f \in (Q_i \setminus pT)R[u]$ we can write $f = qg$ where $q \in Q_i \setminus pT$ and $g \in R[u].$ Find an $n$ such that $p^n f \in R$ and an $n'$ such that $p^{n'}g \in R.$ Let $m= \max \{n, n'\}$ so that $p^m f \in R$ and $p^m g \in R.$ But we see that $p^m f = q p^m g,$ and so $p^m f \in ((Q_i \setminus pT)R) \cap R,$ and since $R$ is a $p$ca subring we know $p^m f=0.$  Since $ p$ is not a zero divisor we have $f=0$ and so $f \in ((Q_i \setminus pT)R[u])\cap R[u]=\{0\}.$

Now we check that we maintain this property when localizing at $R[u] \cap \m$ which will give us that  $((Q_i \setminus pT)S)\cap S = \{0\}$ for any  $i \in \{1, 2, \ldots, k\}$.  Suppose we have an element $ s \in ((Q_i \setminus pT)S)\cap S.$ We can then write $s = f/g = q f'/g'$ with $f , g, f', g' \in R[u]$ with $g$ and $g'$ units in $T$ and $q \in Q_i \setminus pT$. Then we have $fg' = q f' g$ and so clearly $f g' \in ((Q_i \setminus pT)R[u])\cap R[u] $ and so $fg' =0.$ Since $g'$ is a unit, we have $f = 0$ and thus $s=0.$ We have now shown $((Q_i \setminus pT)S)\cap S = \{0\}.$

Finally, let $P \in \Ass T$ and let $f\in
P\cap R[u]$.  Choose an $n$ such that $p^n f \in R$. Then $p^n f \in R \cap P$ and so $p^n f =0$ since $R$ is a $p$ca subring.  Since $p$
is not a zerodivisor, $f=0$ and so we have that $P\cap S=0$. We have now verified all the conditions necessary to show that $S$ is a $p$ca subring of $T$.
\end{proof}

In the constructions in the sequel, we will often need to take unions of $p$ca subrings at intermediate steps. The purpose of Lemma \ref{unions} is to avoid repeating the arguments checking that the union is still a $p$ca subring.

\begin{lemma}
\label{unions}
Suppose we have $(T, \m),$ $C,$ and $p$ as in the hypotheses of Lemma \ref{Claim1}.
Let $\Omega$ be a well ordered set and let $\{R_\alpha | \alpha \in \Omega \}$ be a set of $p$ca subrings indexed by $\Omega$ with the property $R_\alpha \subseteq R_\beta$ for all $\alpha$ and $\beta$ such that $\alpha<\beta.$ Let $S =\bigcup_{\alpha \in \Omega} R_\alpha.$ Then $S \cap P =(0)$ for all associated primes $P$ of $T$ and  $((Q_i \setminus pT) S) \cap S= \{ 0 \}$ for each $i \in \{1,2, \ldots, k\}$.  Furthermore if $\Gamma(R_\alpha) \leq \lambda$ for all $\alpha \in \Omega$ we have $\Gamma(S) \leq \lambda \Gamma(\Omega)$ and so if $\Gamma(\Omega) \leq \lambda$ and $\Gamma(R_\alpha) = \lambda$ for some $\alpha$ we have $\Gamma(S) = \lambda.$
\end{lemma}
\begin{proof}
No further explanation is necessary for the cardinality conditions. Clearly $S \cap P = (0)$  for all $P \in \mbox{Ass}T$ because the $R_\alpha $ are $p$ca  subrings and so none contains a nonzero element of any associated prime ideal of $T.$ Finally, suppose we have $((Q_i \setminus pT)S )\cap S\not= \{0\}$ for some $Q_i \in C.$ Then for some $r, \ r' \in S$ with $r \neq 0$ and for some $q \in Q_i \setminus pT$ we have $r = q r'$.  If we choose $\alpha$ such that $r, \ r' \in R_\alpha$, we see $((Q_i \setminus pT)R_\alpha )\cap R_\alpha \not= \{0\}$, contradicting the hypothesis that $R_\alpha$ is a $p$ca subring.
\end{proof}

\begin{definition}
Let $\Omega$ be a well ordered set and $\alpha \in \Omega$.  We define
$\gamma(\alpha)=\sup\{\beta\in\Omega\suchthat\beta<\alpha\}$.
\end{definition}

\begin{lemma}
\label{Claim2}
Suppose we have $(T, \m),$ $C,$ and $p$ as in the hypotheses of Lemma \ref{Claim1}. Given $(R,R \cap \m)$ a $p$ca subring of $(T, \m),$   there exists a $p$ca subring
$S$ of $T$ with $\Gamma(S)=\Gamma(R)$ such that $R\subseteq S \subseteq T$ and $pT \cap
R\subseteq pS$.
\end{lemma}

\begin{proof}

Let $\Omega=pT \cap R$.  Clearly $\Gamma(\Omega ) \leq \Gamma(R)$.  Well order $\Omega$ and let $0$ denote the first element.  Let
$R_0=R$ and $\alpha \in \Omega$.  Using induction, assume that a $p$ca subring $R_\beta$ with $\Gamma(R_\beta) = \Gamma(R)$ has been defined for every $\beta <\alpha$ so that $\delta  \in pR_{\beta} $ for all $\delta < \beta.$

If $\gamma(\alpha)<\alpha$, then
construct $R_\alpha$ from $R_{\gamma(\alpha)}$ using Lemma \ref{Claim1}
with $c = \gamma(\alpha)$.  By
construction, $R_\alpha$ is a $p$ca subring of $T$ and $\Gamma(R_\alpha )= \Gamma(R_{\gamma(\alpha)})=\Gamma(R)$. Since $\gamma(\alpha) \in p R_\alpha$ and $R_{\gamma( \alpha )} \subseteq R_\alpha,$ using the induction hypothesis, we see that $\delta \in p R_\alpha$ for all $\delta < \alpha.$

Otherwise,
$\gamma(\alpha)=\alpha$ so define $R_\alpha=\bigcup_{\beta<\alpha}
R_\beta$. Then $R_\alpha$ is a union indexed over a segment of $\Omega$ (which can have cardinality at most $|R|$) of $p$ca subrings of cardinality at most $\Gamma(R)$ and so by Lemma \ref{unions} we know $R_\alpha$ is a $p$ca subring and $\Gamma(R_\alpha)=\Gamma(R).$   Since $\sup\{\beta \in \Omega \suchthat \beta < \alpha\} = \alpha$ in this case, for any $\delta < \alpha$ we can choose a $\beta$ such that $\delta < \beta < \alpha,$ and since by the induction hypothesis $\delta \in p R_\beta$, we know $\delta \in p R_\alpha.$

First, suppose $\Omega$ has no maximal element. Let $S= \bigcup_{\alpha \in \Omega} R_\alpha$.
 $S$ is a union of $p$ca subrings $R_\alpha$ such that $\Gamma(R_\alpha)= \Gamma (R)$ for all $\alpha$ indexed by a set of cardinality at  most $|R|$ so by Lemma \ref{unions} we know $S$ is a $p$ca subring with $\Gamma(S) =\Gamma(R).$  Additionally, if $r\in pT\cap R$ then
$r=\gamma(\alpha)$ for some $\alpha$ in $\Omega$ with
$\gamma(\alpha)<\alpha$.  Thus $r\in pR_\alpha \subseteq pS$, so
$pT\cap R\subseteq pS$, and we see that $S$ is our desired subring.

Otherwise, let $d$ denote the maximal element of $\Omega.$ Construct $S$ from $R_d$ using Lemma \ref{Claim1} with $c = d.$  By
construction, $S$ is a $p$ca subring of $T$ with $R_d \subseteq S,$ $d \in pS$ and $\Gamma(S)= \Gamma(R_d)=\Gamma(R)$. Finally, for every $r \in pT \cap R= \Omega,$ either $r <d$, in which case we know by induction $r \in pR_d$ and so $r \in pS$, or $r=d$, in which case $r \in pS$ by our construction of $S$. So $pT \cap R \subseteq pS$, so we see that $S$ is our desired subring.

\end{proof}

\begin{lemma}
\label{fixing}
Suppose we have $(T, \m),$ $C,$ and $p$ as in the hypotheses of Lemma \ref{Claim1}. Let $(R,R \cap \m)$ be a $p$ca subring of $(T, \m)$.  Then there exists a $p$ca subring $S$ of
$T$ with $\Gamma(S)= \Gamma(R)$ such that $R\subseteq S \subseteq T$ and $pT \cap S=pS$.
\end{lemma}

\begin{proof}
Let $R_0=R$.  We define $R_i$ by induction.  Assuming $R_{i-1}$ has been defined so that it is a $p$ca subring and $\Gamma(R_{i-1})=\Gamma(R)$, we use
Lemma
\ref{Claim2} to find a $p$ca subring $R_i$ with $pT\cap R_{i-1}
\subseteq
pR_i$ and $\Gamma(R_i) = \Gamma(R_{i-1})=\Gamma(R).$  Let $S=\bigcup_{i=1}^\infty R_i$.  By Lemma \ref{unions} we know $S$ is a $p$ca subring with $\Gamma(S) =\Gamma(R).$   Further, if $c\in pT\cap S,$ there is an $n\in \N$ such that $c\in pT\cap R_n\subseteq pR_{n+1}\subseteq pS$.  Therefore $pT\cap S \subseteq pS$, so $pT\cap S= pS$.
\end{proof}




The following is Proposition 1 from \cite{Heit}.  It helps us to
ensure that the final ring we create has $T$ as its completion.

\begin{proposition}
\label{prop}
\begin{normalfont}\textbf{(Heitmann,
\cite{Heit})}\end{normalfont}
If $(R,\m\cap R)$ is a quasilocal
subring of a complete local ring $(T,\m )$, the map $R\rightarrow
T/\m^2$ is onto and $IT\cap R=I$ for every finitely generated ideal $I$
of $R$, then $R$ is Noetherian and the natural homomorphism
$\widehat{R}\longrightarrow T$ is an isomorphism.
\end{proposition}

We will construct $A$ so that the map
$A\rightarrow T/\m^2$ is onto. To do this, we will need
Lemma \ref{New3.5}, which lets us adjoin an element of a coset of $T/J$
to a $p$ca subring $R$ where $J$ is an ideal of $T$ such that $J \nsubseteq Q_i$
for every $i \in \{1,2, \ldots,k\}$
 to get a new $p$ca subring.
With $J = \m^2,$ we will get that $A \rightarrow T/\m^2$ is onto as desired. Note that Lemma \ref{New3.5} is similar in purpose to Lemma 3.9 of \cite{SM03}.

\begin{lemma}
\label{New3.5} Let $(T,\m)$ be a complete local ring with $\dim T \geq 1$ and suppose we have a finite, pairwise incomparable set of nonmaximal ideals $C=\{Q_1, Q_2, \ldots, Q_k\} \subseteq \Spec T$.  Let $p \in \bigcap_{i=1}^{k} Q_i$ be a nonzero regular element of $T$ such that for every $P \in { \Ass(T/pT)}$  we have $P \subseteq \bigcup_{i=1}^k Q_i.$

Let $(R,R \cap \m)$ be a $p$ca subring of $T$ such that $pT \cap R =pR$
and let $u+J \in T/J$ where $J$ is an ideal of $T$ with $J\not\subseteq Q_i$ for all $i \in \{1,2, \ldots ,k\}$.
Then there exists an infinite $p$ca subring $S$ of $T$ meeting the following conditions:
\begin{enumerate}
\item $R \subseteq S \subseteq T$
\item $\Gamma(S)=\Gamma(R)$
\item $u+J$ is in the image of the map $S \to T/J$
\item If $u \in J$, then $S \cap J \nsubseteq Q_i$ for each $i \in \{1,2, \ldots ,k\}$
\item $pT \cap S =pS$.
\end{enumerate}
\end{lemma}

\begin{proof}

For each $P \in \Ass T$, let $D_{(P)}$ be a full set of coset
representatives of the cosets $t+P$ that make $(u+t)+P$ algebraic
over $R$. For each $i \in \{1,2, \ldots, k\}$, let $D_{(Q_i)}$ be a full set of coset representatives of
the cosets $t+Q_i \in T/Q_i$ with $t\in T$ that make $(u+t)+Q_i$ algebraic over
$R/R\cap Q_i$ (note that there is no conflict of terminology because $C \cap \Ass T = \emptyset$ since every prime ideal in $C$ contains the regular element $p$). Let $G$ be the set $C \bigcup \{ P_1,P_2, ... , P_k\}$ where
$\{P_1, P_2,\ldots, P_k\} =\Ass T$ and let $D:=\bigcup_{P \in G
}D_{(P)}$. By Lemma 2.3 of \cite{SLGFF} we know since $\dim T \geq 1$ that $|T| \geq |\mathbb{R}|.$ Thus, because $\Gamma(R) < |T| $ we have $|R| < |T|$ and so $|D_{(P)}|< |T|$ for all $P \in G$, and thus we
have that $|D|<|T|$.

We can now employ Lemma 2.4 of \cite{SLGFF} with $I=J$ to find an $x\in J$
such that $x\notin \bigcup \{ r+P\suchthat r\in D, P \in G \}$ since
the set $ C \bigcup \{P_1, \ldots, P_k\}$ is finite. We claim that $S' =
R[u+x]_{(R[u+x] \cap \m)}$ is an infinite $p$ca subring. It's clear that $S'$ is infinite and $\Gamma(S') =\Gamma(R)$.
Further, note that since $(u+x) + P$ is transcendental over $R$ for all $P \in \Ass T$  we know if $f = r_n (u+x)^n + \cdots + r_1 (u+x) + r_0 \in R[u+x] \cap P$ for some $P \in \Ass T$ then $r_i =0$ for every $i$ and so $f=0.$ We thus have $ S'\cap P = (0)$ for every $P \in \Ass T.$

Finally, we claim that for each $i \in \{1,2, \ldots, k\}$,  $((Q_i \setminus pT )S')\cap S' = \{0\}.$ First suppose we have for some $i$ an $f \in ((Q_i \setminus pT) R[u+x]) \cap R[u+x]$ with $f \neq 0.$ Then we have $f = r_n (u+x)^n + \cdots + r_1 (u+x) + r_0= q (s_{n'}(u+x)^{n'} + \cdots +s_1 (u+x) +s_0)$ for some $q \in Q_i \setminus pT$ and some $r_0, r_1, \ldots, r_n, s_0, s_1, \ldots, s_{n'} \in R$ with $r_k \neq 0$ for some $1 \leq k  \leq n.$ Let  $m$ be the largest integer such that $r_i \in (pT)^m$ for all $1\leq i \leq n$ and let $m'$ be the largest integer such that $s_j \in (pT)^{m'}$ for all $1 \leq j \leq n' .$ Then since $pT \cap R = pR$ we have $(pT)^m \cap R = p^m R $ (and similarly for $m'$) and we can write $f = p^m (r_n' (u+x)^n + \cdots + r_1' (u+x) + r_0') =q p^{m'}(s_{n'}'(u+x)^{n'} + \cdots + s_1' (u+x) +s_0')$ for some $r_0', r_1, \ldots, r_n', s_0', s_1', \ldots, s_{n'}' \in R.$

By the maximality of $m$ and $m'$ we know there is an $l$ such that $r_l' \notin pT$ and a $j$ such that $s_j' \notin pT.$ Since $((Q \setminus pT)R) \cap R =(0)$ for all $Q \in C$ we know $Q \cap R \subseteq pT$ and thus $r_l', s_j' \notin Q\cap R$ for all $Q \in C.$ Since $(u+x)+ Q $ is transcendental over $R/R \cap Q$ for all $Q \in C$ we therefore know that $r_n' (u+x)^n + \cdots + r_1' (u+x) + r_0' \notin \bigcup_{i=1}^k Q_i$ and $s_{n'}'(u+x)^{n'} + \cdots +s_1' (u+x) +s_0' \notin \bigcup_{i=1}^k Q_i .$ Now suppose that $m \leq m'.$ Since $p$ is not a zero divisor we may cancel it on both sides of our equation to get $r_n' (u+x)^n + \cdots + r_1' (u+x) + r_0'=q p^{m'-m}(s_{n'}'(u+x)^{n'} + \cdots+ s_1' (u+x) +s_0').$ But the left hand side is not in $\bigcup_{i=1}^k Q_i$ while the right hand side is clearly in $Q_i,$ which is a contradiction. On the other hand suppose $ m > m'.$ Then canceling we have  $p^{m-m'} (r_n' (u+x)^n + \cdots + r_1' (u+x) + r_0') =q( s_{n'}'(u+x)^{n'} + \cdots+ s_1' (u+x) +s_0').$ The left hand side is clearly in $pT$ but since $s_{n'}'(u
 +x)^{n'} + \cdots +s_1' (u+x) +s_0'$ is not in $\bigcup_{i=1}^k Q_i$ it is not in any associated prime of $pT$ and so is not a zero divisor of $T/ pT.$ Since $q \notin pT$ we have that the right hand side is not in $pT,$ which is a contradiction. Thus we have $((Q_i \setminus pT) R[u+x] )\cap R[u+x]=\{0\}.$  By the same trivial checking performed in the proof of Lemma \ref{Claim1} we know that localizing preserves this property and so $((Q_i \setminus pT )S')\cap S' = \{0\}$ for all $Q_i \in C.$ We have now shown that $S'$ is a $p$ca subring of $T.$

 We now employ Lemma \ref{fixing} to find a $p$ca subring $S$ with $S' \subseteq S \subseteq T$ and $\Gamma(S)=\Gamma(S')=\Gamma(R)$ such that $pT \cap S = pS.$ Since $S' \subseteq S$, the image of $S$ in $T/J$ contains $ u+x +J =u+J.$ Furthermore, if $u \in J$ then $u+x \in J \cap S$ but since $(u+x)+Q_i$ is transcendental over $R/R\cap Q_i$ for each $i \in \{1,2, \ldots, k\}$, we have $u+x \notin Q_i$ so $J \cap S \nsubseteq Q_i$ for all $i$.
\end{proof}

The following two lemmas, which are similar to Lemmas 3.10 and 3.11 of \cite{SM03}, allow us to construct A such that  $IT\cap A=I$ for every
finitely generated ideal $I$ of $A$. This is one of the conditions
from Proposition \ref{prop} needed to show that $\widehat{A}=T$.

\begin{lemma}
\label{New3.6} Suppose we have $(T, \m),$ $C,$ and $p$ as in the hypotheses of Lemma \ref{New3.5}. Let $(R,R \cap \m)$ be a $p$ca subring of $(T,\m)$ such that $pT \cap R =pR,$  let $I$ be a finitely
generated ideal of $R,$  and let $c\in IT\cap R$. Then there exists a
$p$ca subring $S$ of $T$ meeting the following conditions:

\begin{enumerate}
\item $R \subseteq S \subseteq T$
\item $\Gamma(S)=\Gamma(R)$
\item $c\in IS$
\item  $pT \cap S = p S$.
\end{enumerate}
\end{lemma}

\begin{proof}
We first show that there exists a $p$ca subring $S'$ of $T$ satisfying the first three conditions.  Induct on the number of generators of $I$. Suppose $I=aR$. If $a=0$, then $c=0$ so $S'=R$ is the desired $p$ca subring. If $a \neq
0$, then $c=au$ for some $u\in T$. We claim that $S'=R[u]_{(R[u]\cap
\m)}$ is the desired subring. First note that clearly $\Gamma (S' ) =\Gamma(R) < \Gamma (T)$. Let $P \in \Ass T$ and suppose $f \in P$. Then $f=r_nu^n+\cdots+r_1u+r_0\in P$,
and $a^nf=r_nc^n+\cdots +r_1ca^{n-1}+r_0a^n\in P\cap R=(0)$. Since
$a \in R$ and $R$ contains no zero divisors of $T$, $f=0$ and so
$S' \cap P = 0 $. Now suppose $f \in ((Q_i \setminus pT ) R[u] )\cap R[u]$ for some $i \in \{1, 2, \ldots ,k\}.$   Then $f = q g$ where $q \in Q_i \setminus pT$ and  $ g \in R[u].$ Since $c=au \in R$, from the argument above we know we have an $m$ such that $a^m f \in R$ and $a^m g\in R.$ Thus we have $a^m f \in ((Q_i\setminus pT)R) \cap R$ and since $R$ contains no zero divisors of $T,$ we know $f =0.$ Therefore $((Q_i \setminus pT ) R[u]) \cap R[u]=\{0\}$ for all $i$. A trivial checking as in the proof of Lemma \ref{Claim1} now verifies that $((Q_i \setminus pT ) S')\cap S'=\{0\}.$

Now let $I$ be an ideal of $R$ that is generated by $m>1$ elements,
and assume that the lemma holds for all ideals with $m-1$ generators.
Let $I=(y_1,\ldots, y_m)R$. Since $c \in IT$ we can choose $t_1, t_2, \ldots, t_m \in T$ such that $c=y_1t_1+y_2t_2+\cdots+y_mt_m.$

First suppose that $y_j \not\in pT \cap R = pR$ for some $j = 1,2, \ldots ,m$.  Without loss of generality, reorder the $y_i$'s so that $y_2 \not\in pT \cap R$.
Our goal is now to find a $t \in T$ such that we may adjoin $t_1+y_2t$ to our subring $R$ without disturbing the $p$ca properties. First note that if $(t_1+y_2t)+Q_i=(t_1+y_2t')+Q_i$
for any $i$, then we have that $y_2(t-t')\in Q_i$. However by the assumption that $y_2 \notin pR$ and the fact that $ Q_i \cap R=pT \cap R =pR,$ we know that $y_2\notin Q_i.$ Since
$Q_i$ is prime, we must have $(t-t')\in Q_i$, thus
$t+Q_i=t'+Q_i$. Therefore, if $t+Q_i\neq
t'+Q_i$, then $(t_1+y_2t)+Q_i\neq (t_1+y_2t')+Q_i.$

For each $i$ let $D_{(Q_i)}$ be a
full set of coset representatives of the cosets $t+Q_i$ that make
$t_1+ y_2 t+Q_i$ algebraic over $R/R\cap Q_i$. Also for any $P \in \Ass T$ let
$D_{(P)}$ be a full set of coset representatives of the cosets $t+P$
that make $t_1+ y_2 t +P$ algebraic over $R$  (again, note that there is no conflict of terminology because $C \cap \Ass T = \emptyset$).
Let $C'$ be the set $\{ Q_1, Q_2, \ldots, Q_k, P_1, P_2, \ldots , P_{k'}\}$ where $\{P_1, P_2,\ldots, P_{k'}\} =\Ass T.$
Let $D=\bigcup_{P \in C} D_{(P)}$. Note that $\m \notin C$ since $Q_i \neq \m$ for all $i$ by assumption and $\m \notin \Ass T$ from the fact that $p \in \m$ is regular. Using the fact from the previous paragraph that $(t_1+y_2t)+Q_i\neq (t_1+y_2t')+Q_i$ whenever $t+Q_i\neq
t'+Q_i,$ it can be easily checked that $|D| < |T|$ and thus we  use Lemma 2.4 of \cite{SLGFF} with $I=T$ to find an
element $t\in T$ such that $t\notin \bigcup \{r+P\suchthat r\in D, P
\in C\}$. Thus, letting $x =t_1+ y_2 t$  we have that $x+Q_i$ is transcendental over $R/R\cap
Q_i$ for all $i$ and $x + P$ is transcendental over $R$ for all $P \in \Ass T $. We now know that $R':=R[x]_{(R[x]\cap \m)}$ is a $p$ca subring of $T$ by the argument in the proof of Lemma \ref{New3.5}.

 We now both add and subtract $y_1 y_2 t$ to see that
$c=y_1t_1+y_1y_2t-y_1y_2t+y_2t_2+\cdots+y_mt_m=y_1 x+y_2(t_2-y_1t)+y_3t_3+\cdots
+y_mt_m.$
Let $J=(y_2,\ldots,y_m)R'$ and $c^*=c-y_1x$. Then $c^*\in JT\cap
R'$ and so we use the induction assumption to find a $p$ca subring $S'$
of $T$ with $\Gamma(S') = \Gamma (R) $ such that $R'\subseteq S\subseteq T$ and $c^* \in JS$. Then
$c=y_1x+c^*\in IS'$, and $S'$ is our desired
$p$ca subring.

Now suppose that $y_j\in pT \cap R$ for all $j$. Then let $k$ be the largest integer such that $y_j \in (pT)^k \cap R $ for all $j$. Since $pT \cap R = pR$ we know $(pT)^k \cap R =p^k R$ and we can write $c=p^k(y_1't_1+\cdots +y_m' t_m)$ for some $y_1', y_2', \ldots, y_m' \in R$ such that $y_i' \notin pT$ for some $i.$ Now let $I' = (y_1', \ldots, y_m')R$ so that we have $c' :=y_1't_1+\cdots +y_m' t_m  \in I' T.$ We can now apply the argument above to find a $p$ca subring $S'$ such that $c' \in I' S'$ and so $c' = y_1' s_1+ \cdots +y_m' s_m  $ for some $s_1, \ldots, s_m  \in S'.$ Then we have $c =p^k c' = p^k y_1' s_1' + \cdots+ p^k y_m' s_m = y_1s_1 + \cdots + y_m s_m$ and so $c \in IS'$ showing that $S'$ is our desired $p$ca subring.

Now we apply Lemma \ref{fixing} to find an $p$ca subring $S$ with $R \subseteq S' \subseteq S$ and $\Gamma(S)=\Gamma(S')=\Gamma(R)$ such that $pT \cap S = pS.$ We know $c \in IS$ since $c \in IS'$ and $S' \subseteq S.$ Thus $S$ is a $p$ca subring meeting the conditions stated in the lemma.
\end{proof}

Lemma \ref{New3.7} allows us to create a subring $S$ of $T$ that
satisfies many of the conditions we want to be true for our final
ring $A$.

\begin{lemma}
\label{New3.7}
Suppose we have $(T, \m),$ $C,$ and $p$ as in the hypotheses of Lemma \ref{New3.5}. Let $(R,R \cap \m)$ be a $p$ca subring of $T$ such that $pT \cap R = pR$ and let $J$ be an ideal of $T$ with $J \nsubseteq Q_i$ for all $i \in \{1,2, \ldots, k\}$ and let $u+J \in T/J$.  Then there exists a $p$ca subring $S$
of $T$ such that
\begin{enumerate}
\item $R \subseteq S \subseteq T$
\item $\Gamma(S)=\Gamma(R)$
\item $u+J$ is in the image of the map $S \to T/J$
\item If $u \in J$, then $S \cap J \nsubseteq Q_i$ for each $i \in \{1,2, \ldots ,k\}$
\item For every finitely generated ideal $I$ of $S$, we have $IT \cap S =I$.

\end{enumerate}
\end{lemma}

\begin{proof}
We first apply Lemma \ref{New3.5} to find an
infinite $p$ca subring $R'$ of $T$ satisfying conditions 1, 2, 3, and 4 and such that $pT \cap R' =p R'$. We will now construct the desired $S$ such that $S$ satisfies conditions 2, and 5 and $R' \subseteq S \subseteq T$ which will ensure that the first, third, and fourth
conditions of the lemma hold true.  Let $\Omega = \{ (I,c) | I$ is a finitely
generated ideal of $R'$ and $c \in IT \cap R' \}$.  Letting $I = R'$,
we see that $| \Omega | \geq | R' |$.  Since $R'$ is infinite, the
number of finitely generated ideals of $R'$ is $| R' |$, and therefore
$| R' | \geq | \Omega |$, giving us the equality
$| R'| = | \Omega |$ and thus $\Gamma( \Omega) = \Gamma (R).$ Well order $\Omega$ so that it does not have a maximal element and
let $0$
denote its first element.  We will now inductively define a family of
$p$ca subrings of $T$, one for each element of $\Omega$.  Let $R_0 =
R'$, and let
$\alpha \in \Omega$.  Assume that $R_\beta$ has been defined for all
$\beta< \alpha$ and that $pT \cap R_\beta = p R_\beta$ and $\Gamma(R_\beta ) = \Gamma(R)$ hold for all $\beta < \alpha.$  If $\gamma(\alpha) < \alpha$ and $\gamma(\alpha) = (I,c)$, then
define $R_\alpha$ to be the $p$ca subring obtained from Lemma
\ref{New3.6}.  Note that clearly $pT \cap R_\alpha = p R_\alpha$ and $\Gamma(R_\alpha )= \Gamma(R_{\gamma(\alpha)}) = \Gamma(R).$
If on the other hand $\gamma(\alpha) = \alpha$, define $R_\alpha = \bigcup_{\beta <
\alpha} R_\beta$.  By Lemma \ref{unions} $R_\alpha$ is a $p$ca subring with $\Gamma(R_\alpha)=\Gamma(R).$ Furthermore, if $t \in pT \cap R_\alpha$ then $t \in R_\beta$ for some $\beta<\alpha$ and so $t \in pT \cap R_\beta= pR_\beta \subseteq  pR_\alpha.$ Thus $pT \cap R_\alpha = pR_\alpha.$

Now let $R_1 = \bigcup_{\alpha \in \Omega} R_\alpha$. We see from Lemma \ref{unions} that $R_1$ is a $p$ca subring and $\Gamma(R_1) = \Gamma(R_0) =\Gamma(R).$ Also, since we know by induction that $pT \cap R_\alpha =p R_\alpha$ for all $\alpha \in \Omega$ we see by the same argument made at the end of the last paragraph that $pT \cap R_1 = p R_1.$
Furthermore, notice that if $I$ is a
finitely generated ideal of $R_0$ and $c \in IT \cap R_0$, then
$(I,c) = \gamma(\alpha)$ for some $\alpha \in \Omega$ with
$\gamma(\alpha) < \alpha$.  It follows from the construction that $c
\in IR_\alpha \subseteq IR_1$.  Thus $IT \cap R_0
\subseteq IR_1$ for every finitely generated ideal $I$ of $R_0$.

Following this same pattern, build a $p$ca subring $R_2$ of $T$ with $\Gamma (R_2) =\Gamma(R_1)=\Gamma(R)$  and $pT \cap R_2 = pR_2$ such
that $R_1 \subseteq R_2 \subseteq T$ and $IT \cap R_1
\subseteq IR_2$ for every finitely generated ideal $I$ of $R_1$.
Continue by induction, forming a
chain $R_0 \subseteq R_1 \subseteq R_2 \subseteq \cdots$ of
$p$ca subrings of
$T$
such that $IT \cap R_n \subseteq IR_{n+1}$ for every finitely
generated
ideal
$I$ of $R_n$ and $| R_i | = | R_0 |$ for all $i$.

We now claim that $S = \bigcup_{i =1}^\infty R_i$ is the desired
$p$ca subring.
To see this, first note $R \subseteq S \subseteq T$ and that we know from Lemma \ref{unions} that $S$ is indeed a $p$ca subring and $\Gamma(S) = \Gamma (R)$.  Now set $I = (y_1, y_2, \ldots , y_k)S$
and let $c \in IT \cap S$.  Then there exists an $N \in \N$ such that
$c, y_1, \ldots , y_k \in R_N$.
Thus $c \in (y_1, \ldots , y_k)T \cap R_N \subseteq (y_1, \ldots ,
y_k)R_{N+1}
\subseteq IS$.  From this it follows that $IT \cap S = I$, so the
fifth condition of the statement of the lemma holds.
\end{proof}

In Lemma \ref{New3.8} we construct a domain $A$ that has the desired
completion and the formal fiber of $pA$ is semilocal with maximal ideals
the elements of $C$.

\begin{lemma}
\label{New3.8}
Suppose we have $(T, \m),$ $C,$ and $p$ as in the hypotheses of Lemma \ref{New3.5}. Let $\Pi$ denote the prime subring of $T.$  Suppose $((Q_i\setminus pT)\Pi[p])\cap \Pi [p] =\{0\}$  for all $i$ and that $P \cap \Pi[p]= (0) $ for all $ P \in \Ass T.$ Then there exists a local domain $A\subseteq T$ such that
\begin{enumerate}
\item $p \in A$
\item $\widehat{A} = T$
\item $pA$ is a prime ideal in $A$ and and has a semilocal formal fiber with maximal ideal the elements of $C$
\item If $J$ is an ideal of $T$ satisfying $J \not\subseteq Q_i$ for all $i \in \{1,2, \ldots ,k\}$
then the map $A \rightarrow T/J$ is onto and $J \cap A \nsubseteq Q_i$ for all $i \in \{1,2, \ldots ,k\}.$
\end{enumerate}
\end{lemma}
\begin{proof}
Let $\Omega = \{ u+J \in T/J | J$ is an ideal of $T$ with $J
\nsubseteq Q_i$ for all $i \in \{1, 2, \ldots ,k\}$ and for each $\alpha \in \Omega$ define $\Omega_\alpha :=\{\beta \in \Omega | \beta \leq \alpha\}$.  Since $T$ is infinite and Noetherian, $|\{ J$ is an ideal
of $T$ with $J \nsubseteq Q\}|\leq |T|$.  Also, if $J$ is an ideal of
$T$, then $|T/J|\leq |T|$.  It follows that $|\Omega | \leq |T|$.
Well order $\Omega$ so that each element has fewer than $|
\Omega |$ predecessors.  Let $0$ denote the first element of
$\Omega$.  Define $R'_0$ to be $\Pi[p]$ localized at $\Pi[p] \cap \m$.  We know $\Gamma(R'_0)= \aleph_0$ and since $\dim T \geq 1$ by Lemma 2.3 of \cite{SLGFF} we know that $|T| \geq |\mathbb{R}|$ and thus $\Gamma(R'_0) < |T|.$ Now we can use the same checking argument that is in the proof of Lemma  \ref{Claim1} to see that $R_0'$ is a $p$ca subring of $T.$ We now apply Lemma \ref{fixing} to find
a $p$ca subring $R_0''$ with  $R_0' \subseteq R_0''$ such that $pT \cap R_0'' = p R_0''$ and $\Gamma(R_0'')= \Gamma(R_0') = \aleph_0.$ Next apply Lemma \ref{New3.7} with $J=T$ to find a $p$ca subring $R_0$ with $R_0'' \subseteq R_0$ such that $I T \cap R_0 = I$ for every finitely generated ideal $I$ of $R_0$ and $\Gamma(R_0)= \Gamma(R_0'') =  \aleph_0.$

Starting with $R_0,$ recursively define a family of $p$ca subrings as follows.  Let $\alpha \in \Omega$ and assume that
$R_\beta$ has already been defined to be a $p$ca subring for all $\beta < \alpha$ with $IT \cap R_\beta = I R_\beta$  for every finitely generated ideal $I$ of $R_\beta$ and $\Gamma(R_\beta) \leq \Gamma(\Omega_\beta)$ (note that this condition holds for $R_0$ since $\Gamma(R_0)= \Gamma(\Omega_0) = \aleph_0$).   Then $\gamma(\alpha) = u+J$ for some ideal $J$ of $T$ with $J \nsubseteq Q_i$ for every $i \in \{1,2, \ldots ,k\}$.  If $\gamma(\alpha) < \alpha$, use Lemma
\ref{New3.7} to obtain a  $p$ca subring $R_\alpha$ with $\Gamma(R_\alpha )= \Gamma (R_{\gamma(\alpha)})$ such that
$R_{\gamma(\alpha)} \subseteq R_\alpha \subseteq T$, $u+J$ is in
the image of the map $R_\alpha \to T/J$ and $IT \cap R_\alpha = I$
for every finitely generated ideal $I$ of $R_\alpha$.  Moreover,
this gives us that $R_\alpha \cap J \not\subseteq Q_i$ for every $i \in \{1,2, \ldots ,k\}$ if $u \in J$.  Also, since $\Gamma(R_\alpha) = \Gamma (R_{\gamma(\alpha)})$ and $\Gamma(\Omega_\alpha) = \Gamma(\Omega_{\gamma(\alpha)})$ we have that $\Gamma(R_\alpha)\leq \Gamma(\Omega_\alpha)$.

If $\gamma(\alpha) = \alpha$ ,
define $R_\alpha = \bigcup_{\beta < \alpha} R_\beta$.  Then by Lemma \ref{unions} we see that $R_\alpha$ is a $p$ca subring of $T$. Furthermore we have $\Gamma(R_\beta)\leq \Gamma(\Omega_\beta) \leq \Gamma(\Omega_\alpha)  $ for all $\beta < \alpha$. So by Lemma \ref{unions} we see that $\Gamma(R_\alpha) \leq \Gamma(\Omega_\alpha)$. Now let
 $I = (y_1, \ldots ,
y_k)$ be a finitely generated ideal of $R_\alpha$ and let $c \in IT \cap R_\alpha$.  Then $\{c, y_1, \ldots , y_k\}
\subseteq R_\beta$ for some $\beta < \alpha .$ By the inductive hypothesis, $(y_1,
\ldots, y_k)T \cap R_\beta = (y_1, \ldots , y_k)R_\beta$.  As $c
\in (y_1, \ldots , y_k)T \cap R_\beta$, we have that $c \in (y_1,
\ldots , y_k)R_\beta \subseteq I$.  Hence $IT \cap R_\alpha = I.$

We now know by induction that for each $\alpha \in \Omega,$ $R_\alpha$ is a $p$ca subring with $\Gamma(R_\alpha)\leq \Gamma(\Omega_\alpha)$ and  $I T \cap R_\alpha = I$ for all finitely generated ideals $I$ of $R_\alpha$.
We claim that $A = \bigcup_{\lambda \in \Omega} R_\alpha$ is the
desired domain.

First note that by construction, condition $(4)$ of the lemma is
satisfied.  We now show that the completion of $A$ is $T$.  Note that as $Q_i$ is
nonmaximal in $T$ for all $i$, we have that $\m^2 \nsubseteq Q_i$ for all $i$.
Thus, by the construction, the map $A \to T/\m^2$ is onto.  Furthermore, by an argument identical to the one used to show that $I T \cap R_\alpha = I$ for all finitely generated ideals $I$ of $R_\alpha$ in the case $\gamma(\alpha)= \alpha,$ we know $I'T \cap A = I'$ for all finitely generated ideals $I'$ of $A.$
It follows from Proposition \ref{prop} that $A$ is Noetherian and $\widehat{A}=T$.

Now we show that the formal fiber of $pA$ is semilocal with maximal ideals exactly the ideals of $C.$    We know that if $P \in \Spec T$ with $P \nsubseteq Q_i$ for all $i$ then $P \cap A \nsubseteq Q_i$ for all $i$ and so $P \cap A \neq pA$ which shows that $P$ is not in the formal fiber of $pA.$  Furthermore, since each $R_\alpha$ is $p$ca, by the argument in Lemma \ref{unions} we know that $((Q_i\setminus pT)A)\cap A =\{0\}$ and so in particular  $(Q_i \setminus pT) \cap A = \emptyset$ for all $i$. Thus $Q_i \cap A = pT \cap A = pA$ for each $i$ and so $pA$ is prime and $Q_i$ is in its formal fiber for every $i \in \{1, 2, \ldots ,k\}$. We have now shown the formal fiber of $pA$ is semilocal with maximal ideals exactly the members of $C$.
\end{proof}

Theorem \ref{principal} is our main result.  The previous work in this
section
has been devoted to showing sufficiency of the conditions below, so
the
majority of the proof of Theorem \ref{principal} demonstrates their necessity.

\begin{theorem} \label{principal}
Let $(T,\m)$ be a complete local ring, $\Pi$ the prime subring
of $T$, and $C=\{Q_1, Q_2, \ldots, Q_k\}$  a finite set of non-maximal incomparable prime ideals of $T$.  Let $p \in \bigcap_{i=1}^k Q_k$ with $p \neq 0$. Then there exists a local domain $A \subseteq T$
with $p \in A$  such that $\widehat{A} = T$  and $pA$ is a prime ideal whose formal fiber is semilocal with maximal ideals the elements of $C$ if and only if $P \cap \Pi[p] =(0)$ for all $P \in \Ass T,$ for every $P' \in \Ass(T/pT)$ we have $P' \subseteq Q_i$ for some $i \in \{1,2, \ldots ,k\}$, and $((Q_i \setminus pT) \Pi[p]) \cap \Pi[p]  = \{0\}$ for all $i \in \{1, 2, \ldots , k\}.$

\end{theorem}
\begin{proof}

The condition that $P \cap \Pi[p] =(0)$ for all $P \in \Ass T$ ensures that $p$ is regular. In particular, $Q_1 \notin \Ass T$ and so $Q_1$ has height at least one and so $\dim T \geq 1.$ Because every $P' \in \Ass(T/pT)$ is contained in some $Q_i$ we know $P' \subseteq \bigcup_{i = 1}^{k} Q_i.$  With these observations, Lemma \ref{New3.8} now shows the conditions are sufficient. We must prove they are necessary.

Suppose we have an $A \subseteq T$ with $\widehat {A} =T$ and that $pA$ is prime with a semilocal formal fiber with maximal ideals exactly the elements of $C.$  Since the extension $A \subseteq \widehat{A}=T$ is faithfully flat, any zero divisor of $T$ which is in $A$ must be a zero divisor of $A.$ Since we assume $A$ is a domain, $A$ can contain no such nonzero zero divisor, and in particular, since certainly $\Pi[p] \subseteq A,$ we must have $P \cap \Pi[p] =(0)$ for all $P \in \Ass T$. Furthermore, since the completion of $A/(pT \cap A) =A/pA$ is $T/pT$ we can say that all zero divisors of $T/pT$ (that is, all elements in the image of $\bigcup \Ass T/pT$  under the canonical map $T \to T/pT$) contained in $A/pA$ are zero divisors of $A/pA.$ But $A/pA$ is a domain since $pA$ is prime, thus $A/pA$ cannot contain any nonzero zero divisor of $T/pT$ and so $A$ does not contain any element of $\bigcup \Ass (T/ pT)$ which is not in $pT$. Let $P \in \Ass (T/pT).$
The argument above shows $P \cap A \subseteq pT \cap A = pA$ and since $p \in P$ we also have $pA \subseteq A \cap P $ giving us $P \cap A =  pA.$  Thus $P$ is in the formal fiber of $pA,$ and since we have assumed this formal fiber is semilocal with maximal ideals $\{Q_1, Q_2, \ldots ,Q_k\}$ we know $P \subseteq Q_i$ for some $i.$

Finally, suppose that for some $i$ there is a $f \in ((Q_i \setminus pT) \Pi[p] )\cap \Pi[p] $ with $f \neq 0.$  We know $f = q g$ for some $q \in Q_i \setminus pT$ and some $g \in \Pi[p] \subseteq A$. Since we showed above it is necessary that $P \cap \Pi[p] = (0)$ for all $P \in \Ass T$ and we know $g \neq 0$ since $f \neq 0$, it follows that $g$ is not a zero divisor of $T.$   Now, since $A \subseteq T$ is a faithfully flat extension, we know  $gT \cap A = gA$ and so $gq \in gA $ which implies $q \in A.$ Therefore $Q_i \cap A \nsubseteq pT,$ contradicting the assumption that $Q_i$ is in the formal fiber of $pA.$ Thus it is necessary that $((Q_i \setminus pT) \Pi[p]) \cap \Pi[p]  = \{0\}$ for all $i$.
\end{proof}

\begin{example}

Let $T$ be the complete local ring $\mathbb{R}[[x,y,z,w]]/(x^2-yz)$ and $Q$ be the
non-maximal prime ideal $(x,y,z)$. $T$ is a domain as $(x^2-yz)$ is a prime ideal in
$\mathbb{R}[[x,y,z,w]]$. Note that if $P \in \Ass(T/xT)=\{(x,y), (x,z)\}$ then
$P \subseteq Q$. It is
also the case that $(Q \setminus xT)\Pi[x] \cap \Pi[x] =\{0\}$.    Thus the conditions
of Theorem \ref{principal} are satisfied, and there exists a domain $A$ such that
$\widehat{A}=\mathbb{R}[[x,y,z,w]]/(x^2-yz)$, $xA$ is a prime ideal in $A$, and the
formal fiber of $xA$ is  local with maximal ideal $(x, y, z)$.

\end{example}

\bibliographystyle{plain}
\bibliography{smallbibliography}

\end{document}